\documentclass[12pt]{amsart}
\usepackage{amsfonts}
\usepackage{amssymb}

\theoremstyle{definition}

\theoremstyle{remark}

\begin{document}

\centerline{\bf Comptes rendus de l'Acad\'emie bulgare des Sciences}

\centerline{\it Tome 63, No 1, 2010}

\vspace{0.6in}

\begin{flushright}
{\it MATH\'EMATIQUES
\\ G\'eometrie diff\'eretielle}
\end{flushright}

\vspace{0.2in}

\title{Almost Hermitian manifolds with vanishing Bochner curvature tensor}%
\author{Ognian Kassabov}%
\address{Higher transport School, Sofia, Bulgaria}
\email{okassabov@abv.bg}

\maketitle

\vspace{0.2in}

\centerline{\it (Submitted by Academician P. Popivanov on September 29, 2009)}

\vspace{0.2in}

{\bf Abstract}

\vspace{0.1in}
We prove that if an almost Hermitian manifold of dimension greater than 4 and vanishing 
Bochner curvature tensor is not K\"ahler at a point, then it is flat in
a neighbourhood of this point.



\vspace{0.3in}
{\bf 1. Introduction.} 
In \cite{B}  S. Bochner introduced a curvature tensor as a K\"ahler 
analogue of the Weil conformal curvature tensor. Since then there are many results on the 
Bochner's tensor for K\"ahler manifolds, as well as for some more general classes of
almost Hermitian manifolds, see e.g. \cite{RB}, \cite{VY} and  the references in them. These results show
that the Bochner tensor is in fact a perfect analogue of the Weil tensor, although we don't
know for the moment its geometric meaning.  

For the more general classes  of almost Hermitian manifolds we can use the classical Bochner tensor,
as well as its extension, due to Tricerri and Vanhecke \cite{TV}. For some results concerning
the last tensor see e.g.  \cite{G}, \cite{GM}, \cite{K} and the included references. The main result of \cite{K}
has recently been reproved in \cite{EPS}.   
 
The aim of this paper is to show that the vanishing of the classical Bochner tensor 
provoke very strong restrictions on an almost Hermitian manifold, so it seems more
natural to use in this case the tensor of Tricerri-Vanhecke. Namely we prove the following result:

\vspace{0.1in}
{\bf Theorem.} {\it Let $M$ be an $2n$-dimensional ($n> 2$) almost Hermitian manifold
with vanishing classical Bochner curvature tensor. If $M$ is not K\"ahlerian at a point $p$, 
then $M$ is flat in a neighbourhood of $p$.}

\vspace{0.3in}
{\bf 2. Preliminaries.}
Let $M$ be a $2n$-dimensional almost Hermitian manifold with Rieman\-nian metric $g$ and
almost complex structure $J$. The curvature tensor, the Ricci tensor and the scalar 
curvature are denoted by $R$, $\rho$ and $\tau$, respectively. The Bochner curvature tensor $B$
is defined by
$ B= R - \varphi (Q) $, where
$$
	\begin{array}{rcl}
		\varphi(Q) (x,y,z,u) & = &  g(x,u)Q(y,z) -g(x,z)Q(y,u) \\
                  &	+	&  g(y,z)Q(x,u) -g(y,u)Q(x,z)\\
	                & + &    g(x,Ju)Q(y,Jz) -g(x,Jz)Q(y,Ju) - 2g(x,Jy)Q(z,Ju)  \\
	                & + &    g(y,Jz)Q(x,Ju) -g(y,Ju)Q(x,Jz) - 2g(z,Ju)Q(x,Jy)  \\
	\end{array}
$$
and
$$
   Q = \frac{1}{2(n+2)}\rho -
       \frac{\tau }{8(n+1)(n+2)}g \ .
$$ 

Suppose the Bochner tensor vanishes and let $e_1,...,e_{2n}$ be an orhonormal basis of $T_pM$. 
After adding the identities
$$
	R(x,e_i,e_i,y) = \varphi(Q)(x,e_i,e_i,y)
$$
for $i=1,...,2n$ we obtain
$$
	\rho(x,y) = \frac{1}{2(n+2)}\{(2n+1)\rho(x,y) + 3 \rho(Jx,Jy)\}  \ \ .
$$
Hence the Ricci tensor $\rho$ is hybrid:
$$
	\rho(x,y) =  \rho(Jx,Jy) \ .
$$
So the tensor $Q$ is also symmetric and hybrid:
$$
  Q(x,y) = Q(y,x)  \quad \mbox{and } \quad
  Q(x,y) = Q(Jx,Jy)    \ .  \leqno (2.1)
$$
As a direct consequence of (2.1) $M$ is an $AH_1$-manifold, i.e.
$$
	R(x,y,z,u)=R(x,y,Jz,Ju) \ \ , \leqno (2.2)
$$
see also \cite{VY}. Using (2.2) and the second Bianchi identity 
$$ (\nabla_{x} R)(y,z,u,v) +
   (\nabla_{y} R)(z,x,u,v) +
   (\nabla_{z} R)(x,y,u,v)= 0   \leqno(2.3)
$$
we derive
$$ 
	R(y,z,(\nabla_{x} J)u,Jv) + R(y,z,Ju,(\nabla_{x}J)v)
$$
$$ 
	+ R(z,x,(\nabla_yJ)u,Jv) + R(z,x,Ju,(\nabla_yJ)v) \leqno(2.4)
$$
$$ 
	+ R(x,y,(\nabla_zJ)u,Jv) + R(x,y,Ju,(\nabla_zJ)v)=0 \ \ .
$$
Note also that
$$
	g((\nabla_{x}J)y,y)=0
	\qquad
	g((\nabla_{x}J)y,Jy)=0 \ \ .  
$$

Although the proof of our theorem can be made using real tangent vectors, we
will use the complexification of the tangent spaces. In particular
let $ \{e_{\alpha},\ Je_{\alpha}\},\ \alpha =1,...,n$ be an orthonormal basis
of the tangent space $T_pM$ of $M$ at $p$ and denote:
$$  Z_{\alpha} = e_{\alpha} -iJe_{\alpha} \ , \qquad
    Z_{\bar\alpha} = e_{\alpha} +iJe_{\alpha} \ . $$
Then
$$
	 JZ_{\alpha} = iZ_{\alpha} \ , JZ_{\bar\alpha}=-iZ_{\bar\alpha} \ . \leqno (2.5) 
$$
By (2.1) and (2.5) we find:
$$
	Q(Z_{\alpha},Z_{\beta})=0 \qquad
	Q(Z_{\alpha},Z_{\bar\beta})=Q(Z_{\bar\beta},Z_{\alpha}) \ \ . 
$$

We have also for any $ Z, \ \alpha,\ \beta$
$$
  g((\nabla_ZJ)Z_\alpha,Z_{\bar\beta})=0  \ .
$$

In the following we assume that the basis
$e_{\alpha} ,\ Je_{\alpha},\ \alpha =1,...n$
diagonalize $Q$, i.e.
$$
	Q(e_\alpha,e_\beta)=Q(e_\alpha,Je_\beta)=0 \qquad 
	{\rm for} \qquad \alpha \ne \beta \ .
$$
Then
$$
	Q^1(Z_\alpha)=\mu_\alpha Z_\alpha
	\qquad  Q^1(JZ_\alpha)=\mu_\alpha JZ_\alpha \ \ \ \ \ \alpha = 1,2,...,n 
$$
for some real numbers $\mu_\alpha$, $Q^1$ being the tensor of type (1,1), corresponding to $Q$.

\vspace{0.3in}
{\bf 3. Proof of the theorem.} 
It is sufficient to prove, that if $M$ is non K\"ahlerian in $p$, then it 
is flat at $p$.

First we suppose that there exist  $\alpha,\ \beta$ such that
$$
	g((\nabla_{Z_{\bar\beta}}J)Z_\beta,Z_{\alpha}) \ne 0 \ \ .
$$
We put in (2.4) $X=Z_\alpha,\ Y=Z_{\bar\alpha},\ Z=Z_{\bar\beta},\ U=Z_\alpha,\ V=Z_\beta$ 
and we obtain:
\vspace*{0.05in}
$$
	(5\mu_\alpha+\mu_\beta) g((\nabla_{Z_{\bar\beta}}J)Z_{\beta},Z_{\alpha}) =0  \ \ . \leqno(3.1)
$$
\vspace*{0.1in}
With \ $X=Z_\alpha,\ Y=Z_{\bar\beta},\ Z=Z_{\bar\gamma},\ U=Z_\beta,\ V=Z_\gamma$ \ for
$\gamma \ne \alpha,\ \beta$ (2.4) implies:
$$
	(\mu_\alpha+\mu_\gamma)\ g((\nabla_{Z_{\bar\beta}}J)Z_{\beta},Z_{\alpha})    
	+(\mu_\alpha+\mu_\beta)g((\nabla_{Z_{\bar\gamma}}J)Z_{\gamma},Z_{\alpha}) =0    \ \ .  \leqno(3.2)
$$
Now let in (2.4)  \ $X=Z_{\bar\beta},\ Y=Z_{\gamma},\ Z=Z_{\bar\gamma},\ U=Z_\beta,\ V=Z_\alpha$. The
result is: 
\vspace*{0.05in}  
$$
	(\mu_\alpha+\mu_\beta+2\mu_\gamma) \ g((\nabla_{Z_{\bar\beta}}J)Z_{\beta},Z_{\alpha})    
	-(\mu_\beta+\mu_\gamma) g((\nabla_{Z_{\bar\gamma}}J)Z_{\gamma},Z_{\alpha})=0      \ \ .  \leqno(3.3)
$$

Since \ $	g((\nabla_{Z_{\bar\beta}}J)Z_\beta,Z_{\alpha}) \ne 0 $ \ we can derive from (3.1),(3.2) and (3.3) that 
$M$ is flat at $p$. Indeed, consider the two possibilities:

\vspace{0.1in}
Case I:
$	g((\nabla_{Z_{\bar\gamma}}J)Z_\gamma,Z_{\alpha}) = 0 $. Then (3.1),(3.2) and (3.3) imply
$$
	\left\{ \begin{array}{l}
	  5\mu_\alpha+\mu_\beta=0 \\
	  \mu_\alpha+\mu_\gamma=0 \\
	  \mu_\alpha+\mu_\beta+2\mu_\gamma=0 
	\end{array} \right.
$$
and hence $\mu_\alpha=\mu_\beta=\mu_\gamma=0$.

\vspace{0.1in}
Case II:
$	g((\nabla_{Z_{\bar\gamma}}J)Z_\gamma,Z_{\alpha}) \ne 0 $. By (3.1) we have
$$
	5\mu_\alpha+\mu_\beta=0 \ .
$$
Changing $\beta$ and $\gamma$ we have also
$$
	5\mu_\alpha+\mu_\gamma=0
$$ 
and the last two give $\mu_\beta=\mu_\gamma=-5\mu_\alpha$. Since the system of (3.2) and (3.3)
has a non trivial solution, its determinant is zero:
$$
	(\mu_\alpha+\mu_\gamma)(\mu_\beta+\mu_\gamma)+(\mu_\alpha+\mu_\beta)(\mu_\alpha+\mu_\beta+2\mu_\gamma)=0
$$
Now, using  $\mu_\beta=\mu_\gamma=-5\mu_\alpha$ we derive easily  $\mu_\alpha=\mu_\beta=\mu_\gamma=0$.

Analogously we can see, that
$$
	g((\nabla_{Z_{\beta}}J)Z_\beta,Z_\alpha) \ne 0 
$$
for some \ $\alpha,\ \beta$ \ implies that $M$ is flat at $p$.

Suppose, that there exist $\alpha,\ \beta,\ \gamma$ \ such that
$$
	g((\nabla_{Z_\alpha}J)Z_\beta,Z_\gamma) \ne 0 \ \ .
$$
We put in (2.4) $X=Z_\alpha,\ Y=Z_{\beta},\ Z=Z_{\bar\beta},\ U=Z_\gamma,\ V=Z_\beta$ and we obtain 
\vspace*{0.05in}
$$
	(5\mu_\beta+ \mu_\gamma)\ g((\nabla_{Z_{\alpha}}J)Z_{\beta},Z_{\gamma})    
	-(\mu_\alpha+\mu_\beta) g((\nabla_{Z_{\beta}}J)Z_{\alpha},Z_{\gamma})=0  \ \ .  \leqno(3.4)
$$
\vspace*{0.1in}
Let us change here $\alpha$ and $\beta$:
$$
	(5\mu_\alpha+ \mu_\gamma)\ g((\nabla_{Z_{\beta}}J)Z_{\alpha},Z_{\gamma})    
	-(\mu_\alpha+\mu_\beta) g((\nabla_{Z_{\alpha}}J)Z_{\beta},Z_{\gamma})=0  \ \ .  \leqno(3.5)
$$
Since \ $	g((\nabla_{Z_\alpha}J)Z_\beta,Z_\gamma) \ne 0 $ \ the system (3.4), (3.5) implies
$$
	(5\mu_\alpha+\mu_\gamma)(5\mu_\beta+\mu_\gamma)=(\mu_\alpha+\mu_\beta)^2 \ \ .  
$$
We can change here $\beta$ and $\gamma$: 
$$
	(5\mu_\alpha+\mu_\beta)(5\mu_\gamma+\mu_\beta)=(\mu_\alpha+\mu_\gamma)^2 \ \ .   \leqno(3.6)
$$ 
Consider the following two possibilities:

\vspace{0.1in}
Case I. $g((\nabla_{Z_{\beta}}J)Z_{\alpha},Z_{\gamma})\ne 0$. Then we have also (changing in (3.6) $\alpha$ and $\beta$):
$$
	(5\mu_\beta+\mu_\alpha)(5\mu_\gamma+\mu_\alpha)=(\mu_\beta+\mu_\gamma)^2 \ \ .   
$$ 
The system of the last three equations implies $\mu_\alpha=\mu_\beta=\mu_\gamma=0$.

\vspace{0.in}
Case II. $g((\nabla_{Z_{\beta}}J)Z_{\alpha},Z_{\gamma}) = 0$. Then (3.4), (3.5) imply
$$
	5\mu_\beta+\mu_\gamma = 0 \qquad {\rm and} \qquad \mu_\alpha+\mu_\beta = 0 \ \ ,
$$
so $\mu_\gamma=-5\mu_\beta,\ \mu_\alpha=-\mu_\beta$ and using (3.6) we obtain 
$\mu_\alpha=\mu_\beta=\mu_\gamma=0$.

Consequently   \ $	g((\nabla_{Z_\alpha}J)Z_\beta,Z_\gamma) \ne 0 $ \ always implies  $\mu_\alpha=\mu_\beta=\mu_\gamma=0$.
If $n>3$, we replace in (2.4) $X=Z_\alpha,\ Y=Z_{\delta},\ Z=Z_{\bar\delta},\ U=Z_\beta,\ V=Z_\gamma$ 
for $\delta \ne \alpha,\,\beta,\,\gamma$    and we 
find $\mu_\delta=0$, so $M$ is flat at $p$.

At the end, suppose that there exist $\alpha,\ \beta,\ \gamma$ \ such that
$$
	g((\nabla_{Z_{\bar\alpha}}J)Z_\beta,Z_\gamma) \ne 0 \ \ .
$$
We put in (2.4)  \ $X=Z_{\bar\alpha},\ Y=Z_{\beta},\ Z=Z_{\bar\beta},\ U=Z_\gamma,\ V=Z_\beta$ \ to find
$$
	(5\mu_\beta + \mu_\gamma) g((\nabla_{Z_{\bar\alpha}}J)Z_{\beta},Z_{\gamma}) = 0  \ \ . 
$$
Hence $5\mu_\beta + \mu_\gamma=0$. By the symmetry of $\beta$ and $\gamma$ 
we have also $\mu_\beta + 5\mu_\gamma=0$ and hence $\mu_\beta=\mu_\gamma=0$.

Now we put in the second Bianchi identity (2.3) $X=Z_{\bar\alpha},Y=Z_{\beta},Z=Z_\gamma,U=Z_\beta,V=Z_{\bar\beta}$
and we obtain
$$   
   (\nabla_{Z_{\beta}} Q)(Z_{\bar\alpha},Z_{\gamma}) 
    - 2(\nabla_{Z_{\gamma}} Q)(Z_{\bar\alpha},Z_{\beta})=0  \ \ .
$$
Hence, because of the symmetry in $\beta$ and $\gamma$ we derive easily
\vspace*{0.04in}
$$   
   (\nabla_{Z_{\beta}} Q)(Z_{\bar\alpha},Z_{\gamma}) 
    = (\nabla_{Z_{\gamma}} Q)(Z_{\bar\alpha},Z_{\beta})=0  \ \ . \leqno (3.7)
$$

\vspace*{0.07in}
On the other hand putting in (2.3) $X=Z_{\alpha},Y=Z_{\bar\alpha},Z=Z_\beta,U=Z_\gamma,V=Z_{\bar\alpha}$ 
and using (3.7) we obtain:
$$    
  2(\nabla_{Z_{\bar\alpha}} Q)(Z_{\beta},Z_{\gamma}) 
  + i \mu_\alpha g((\nabla_{Z_{\bar\alpha}} J)Z_{\beta},Z_{\gamma}) =0   \ \ .
$$
Since the first term is symmetric in $\beta,\ \gamma$ and the second is antisymmetric,
they are both zero, so
$   \mu_\alpha  =0   .$

If $n>3$ let $\delta \ne \alpha,\,\beta,\,\gamma$.
We replace in (2.4) $X=Z_{\bar\alpha},\ Y=Z_{\delta},\ Z=Z_{\bar\delta},\ U=Z_\beta,\ V=Z_\gamma$ 
and we derive easily  $\mu_\delta=0$, so $M$ is flat at $p$. 

On the other hand, if 
$$
	g((\nabla_{Z_{\beta}}J)Z_\beta,Z_\alpha) \quad
	g((\nabla_{Z_{\bar\beta}}J)Z_\beta,Z_\alpha) \quad
	g((\nabla_{Z_\alpha}J)Z_\beta,Z_\gamma) \quad
	g((\nabla_{Z_{\bar\alpha}}J)Z_\beta,Z_\gamma) 
$$
are zero for every $\alpha,\,\beta,\,\gamma$, then $M$ is K\"ahlerian at $p$. This proves our theorem.

\end{document}